\documentclass[12pt, amsmath, amsthm, amssymb]{amsart}

\usepackage{xcolor}

\newtheorem{Lemma}{Lemma}[section]
\newtheorem{Theorem}[Lemma]{Theorem}
\newtheorem{Proposition}[Lemma]{Proposition}
\newtheorem{Corollary}[Lemma]{Corollary}

\newtheorem{Assumption}[Lemma]{Assumption}

\newtheorem{remark}[Lemma]{Remark}
\newtheorem{definition}[Lemma]{Definition}
\newtheorem{example}[Lemma]{Example}
\newtheorem{exercise}[Lemma]{Exercise}
\newtheorem{Fact}[Lemma]{Fact}

\def\bt{\begin{Theorem}}
\def\et{\end{Theorem}}
\def\bl{\begin{Lemma}}
\def\el{\end{Lemma}}
\def\bp{\begin{Proposition}}
\def\ep{\end{Proposition}}
\def\bcor{\begin{Corollary}}
\def\ecor{\end{Corollary}}
\def\bpf{\begin{proof}}
\def\epf{\end{proof}}
\def\bex{\begin{exercise}}
\def\eex{\end{exercise}}

\def\brem{\begin{remark}\rm}
\def\erem{\hfill $\lozenge$\end{remark}}

\def\bedef{\begin{definition}}
\def\endef{\end{definition}}

\def\beg{\begin{example}\rm}
\def\eeg{\hfill $\lozenge$\end{example}}

\def\beassum{\begin{Assumption}\rm}
\def\enassum{\hfill $\lozenge$\end{Assumption}}

\def\bef{\begin{Fact}}
\def\eef{\end{Fact}}

\def\bc{\begin{center}}
\def\ec{\end{center}}
\def\noi{\noindent}
\def\vsh{\vskip .5cm}

\def\vsqq{\vskip .125cm}

\def\beq{\begin{equation}}
\def\eeq{\end{equation}}
\def\beqarray{\begin{eqnarray*}}
\def\eeqarray{\end{eqnarray*}}
\def\<{\leftangle}
\def\>{\rightangle}
\def\({\left(}
\def\){\right)}
\def\f{\varphi}

\def\span{{\rm span}}
\def\<{\langle}
\def\>{\rangle}
\def\q{\quad}

\def\a{\alpha}
\def\b{\beta}
\def\g{\gamma}

\def\d{\delta}
\def\h{\hbox}

\def\l{\lambda}

\def\O{\Omega}

\def\m{\mu}

\def\w.r.t.{with respect to}
\def\R{{\mathbb{R}}}
\def\N{{\mathbb{N}}}

\def\C{{\mathbb{C}}}

\def\K{{\mathbb{K}}}

\def\bq{\begin{quote}}
\def\eq{\end{quote}}

\def\dist{{\rm dist}}
\def\bit{\begin{itemize}}
\def\eit{\end{itemize}}
\def\i{\item}
\def\ben{\begin{enumerate}}
\def\een{\end{enumerate}}

\def\H{{\mathcal H}}

\def\X{{\mathcal X}}
\def\Y{{\mathcal Y}}

\def\cb{\color{blue}}

\def\nls{normed linear space }
\def\nlsp{normed linear space}

\begin{document}

%September 21, 2013 (Saturday)
\title[Finite Data Based Learning ]
%{A Finite Data Based Learning in a  \\  Normed Linear Space Setting}
{A Finite Data Based Regularization  in a  \\  Normed Linear Space Setting}

\author{M.Thamban  Nair}
\address{Department of Mathematics, BITS Pilani, K K Birla Goa Campus, Goa 403736, India}
\email{mtnair@goa.bits-pilani.ac.in; mtnair@faculty.iitm.ac.in}
\today
\maketitle

\begin{abstract}  One of the basic  problems in  mathematical learning theory is to identify  a function $f: \O\to \R$ with certain specific properties which fits a given training data  $\{(x_i, \xi_i)\in \O\times \R: i=1, \ldots, n\}$, in the sense that, $f(x_i) =  \xi_i$ for $i=1, \ldots, n$, where $\O$ is a compact subset of $\R^d$ for some $d\in \N$ and $f$ is required to have certain specified characteristics.  We address this problem when $f$ belongs to an arbitrary normed linear space $\X$, and the evaluation maps $f\mapsto f(x_i)$  are replaced by maps of the form  $ f\mapsto \f_i(f)$ on $\X$, where $\f_1, \ldots, \f_n$ are  continuous linear functionals on $\X$.  Using an inner product structure on $\R^n$, we  shall device a method of least-squares for  obtaining an approximate solution for the  above problem  and also identify a subspace of $\X$ in which the least-square solution is unique, which in turn is shown to be equivalent to solving a matrix equation.    In the context when  the  matrix  under consideration is ill-conditioned,  a  regularized equation is devised,  and order optimal error estimates are derived for the exact targeted data $\xi=(\xi_1, \ldots, \xi_n)$ and also when it is noisy, by choosing the regularization parameter appropriately.
\end{abstract} 

\noi 
{\bf Key words:}  Learning theory; reproducing kernel Hilbert space; integral operator;  weighted inner product; least-square solution; regularization, parameter choice. 
\vsqq
\noi 
{\bf AMS Classification:} 45L05, 45Q05, 65J20, 65R30, 65R32

%\tableofcontents 
	
\section{Introduction}

In learning theory, and specifically, in machine learning, the problem that one has to deal with is the following: 
\bq 
Given a training set $\{(x_i, \xi_i)\in \O\times \R: i=1, \ldots, n\}$, where $\O$ is a compact subset of $\R^d$ for some $d\in \N$, look for a function $f: \O\to \R$  with certain specified properties, such that 
$$f(x_i)  = \xi_i, \q i=1, \ldots, n.$$
\eq 
It is essentially an interpolation problem using a specified class of functions. Suppose we are looking for such a function $f$ in the space $C(\O)$, the vector space of all  continuous (real valued) functions  defined on $\O$.
Talking $E_j = \{x_i: i=1, \ldots, n; i\not=j\}$ and 
$$f_j(x) = \frac{d(x, x_j)}{d(x, x_j)+ \dist(x, E_j)},\q x\in \O,$$
where $d(\cdot, \cdot)$ is the usual (euclidean) metric  on $\O$  and   $d(x, S):= \inf_{y\in S}d(x,y)$, we see that each $f_j$ is continuous on $\O$ satisfying $f_j(x_j) =1$ and $f_j(x_i)=0$ for $i\not=j$, and the function $f: \O\to \R$ defined by 
$$f(x) = \sum_{j=1}^n \xi_j f_j(x),\q x\in \O,$$
is continuous and it satisfies $f(x_j) = \xi_j$ for $j=1, \ldots, n$.
In this case, we know that, for each $x\in \O$,  the evaluation map  $\f_x: C(\O)\to \R$, defined by 
$$\f_x(f) = f(x),\q f\in C(\O),$$ 
is a linear functional on $C(\O)$, and it is continuous \w.r.t. the  sup-norm.  
Thus,  we can think of the training set $\{(x_i, \xi_i)\in \O\times \R: i=1, \ldots, n\}$ as 
$$\{(\f_{x_i}, \xi_i)\in \O\times \R: i=1, \ldots, n\}.$$ We note that the above problem cannot be addressed if the function space $C(\O)$ is replaced by  $L^p(\O)$ for some $p$ with $1\leq p\leq \infty$, for in this case evaluation of a function $f\in L^p(\O)$ is not meaningful. However, if we have a   reproducing kernel Hilbert space (RKHS) $\H(\O)$, then  we know that,  for each   $x\in \O$,  $f(x)$ is well-defined for every $f\in \H(\O)$ and the map $f\mapsto f(x)$ is a continuous linear functional (cf. \cite{JH-POA}). In this case,  by Riesz representation theorem,  the evaluation maps can be represented uniquely by functions from $\H(\O)$. Thus, for each $x\in \O$, there exists a unique function $k_x\in \H(\O)$ such that 
$$f(x) = \<f, k_x\>_{\H(\O)}\q\forall\, f\in \H(\O),$$
where $\<\cdot, \cdot\>_{\H(\O)}$ denotes the inner product on $\H(\O)$.   The function $k(\cdot, \cdot)$ defined by $k(x,y) = k_x(y)$ for $(x,y)\in \O\times \O$ is called the {\it reproducing kernel} associated with the space $\H(\O)$. 

The considerations in the  last paragraph motivate us to consider a  learning problem in a more general setting of a normed linear space: 

\bq 
Let  $\X$ be a (real) normed linear  space  and let $\X^*$ denote its dual, that is, the space of all continuous linear functionals on $\X$.  Given a   training set  of the form   $\{(\f_i, \xi_i): i=1, \ldots, n\},$ 
where $\f_i \in  \X^*$ and  $\xi_i\in \R$ for $i=1, \ldots, n$,  we  look for an  $f\in  \X$ such that 
$$\f_i(f)  =  \xi_i\q\h{for}\q  i=1, \ldots, n.$$
\eq 
The above setting can also be considered as generalization of the issues  considered  in the papers  \cite{MTN-SP, Nair-ACM,  Nair-discrete} for obtaining stable approximate solutions for   ill-posed integral  and operator equations,  exploiting  the Hilbert space structure  of certain spaces under consideration.  In order to carry out the analysis in the normed linear settings, in this paper, a novel approach is introduced  which leads to the representation of a {\it least-square solution}, that is,  $f\in \X$ such that 
$$\|\Phi(f)-\xi\|_{\R^n} = \inf_{h\in \X}\|\Phi(h)-\xi\|_{\R^n},$$ 
as a solution of an equation, which is an analogue of the so called {\it normal equation},  of the form 
$$\Phi^*\Phi(f) = \Phi^*(\xi),\q f\in \X,\, \xi\in \R^n,$$ 
where $\Phi: \X\to \R$ is defined by 
$$\Phi(f) = (\f_1(f), \ldots, \f_n(f)),\q f\in \X,$$
and $\Phi^*: \R^n\to \X^*$ is defined by exploiting the inner product structure on $\R^n$, namely, 
$$\Phi^*({\bf u})(f) = \<\Phi(f), {\bf u}\>_{\R^n},\q {\bf u}\in \R^n,\, f\in \X.$$  This formalism enables  us to  identify a unique least-square solution $f^\dagger$ in a certain finite dimensional subspace $\X_0$ of $\X$, which satisfies $\X = N(\Phi)\oplus \X_0$, where $N(\Phi)$ denotes the null space of $\Phi$, and thereby obtain a regularized solution $f_\l$ for $\l>0$  by making use of a  least-square problem and its regularization in the setting of  matrix equations. We shall also derive order optimal error estimates for  the regularized solution when the data is noisy as well as noise-free. 

\brem 
We observe that, for a  given data  $(\Phi, \xi)\in {\mathcal Y}_n\times \R^n$, where ${\mathcal Y}_n:=(\X^*)^n$, the primary goal of the paper is to identify an  $f\in \X$ which  minimizes the function  ${\mathcal E}_{\Phi, \xi}: \X\to \R$, defined by 
$${\mathcal E}_{\Phi, \xi}(f) := \|\Phi(f)- \xi\|_{\R^n}^2,\q f\in \X.$$
This  goal is modest in comparison with an associated   regression problem, where the goal is to identify an  $f\in \X$, which  minimizes the function ${\mathcal E}: \X\to \R$, defined by 
$${\mathcal E}(f) := \int_{{\mathcal Y}_n\times \R}{\mathcal E}_{\Phi, \xi}(f) d\m,$$
where $\m$ is an appropriate probability measure on ${\mathcal Y}_n\times \R^n$.  For a comprehensive analysis of such regression problem, in the setting of  $\X$ being a reproducing kernal Hilbert space  and $n=1$, one may refer to Lin et. al. \cite{lin-et-al}. However,  the problem considered in this paper has relevance to learning theory, and at the same time in includes considerations of  solving ill-posed integral equations and operator equations based on finite data points. Further, the  proposed  regularization procedure, which yields order  optimal error estimates in the setting of a normed linear space. 
\erem

\section{Some examples}   

We consider a few examples to illustrate the general setting of the learning problem.

\beg 
For $1\leq p\leq \infty$, let  $\X$ be a subspace of $L^p(\O)$  and let $\Y = C(\O)$, where $\O$  is a compact subset of $\R^d$ for some $d\in \N$. Let  $K: \X\to \Y$ be   an integral operator with a continuous kernel $k(\cdot, \cdot)$, that is,  
$$(Kf)(x) = \int_\O k(x,y)f(y)dy,\q f\in \X.$$
Suppose we are interested in finding $f\in  \X$  such that 
$$(Kf)(x_i) =  \xi_i,\q i=1, \ldots, n,$$
for a given data set $\{(x_i, \xi_i) \in \O\times \R: i=1, \ldots, n\}$.  Now, defining 
$$\f_i(f) = (Kf)(x_i)\q\h{for}$$
for $f\in \X$ and for $i\in \{1, \ldots, n\}$, we see, using H\"{o}lder's inequality  (cf. Nair \cite{Nair-FA}) that $\f_i$ is a continuous linear functional on $L^p(\O)$ for each $i\in \{1, \ldots, n\}$. 
Thus, the problem here is to identify $f\in L^p(\O)$ such that 
$$\f_i(f) = \xi_i\q\forall\, i\in \{1, \ldots, n\}.$$
In this case, it is to be observed that for each $i\in \{1, \ldots, n\}$, 
$$\f_i(f) = \<f, k_{x_i}\>_{L^2(\O)},$$
where $k_{x_i}(y) = k(x_i, y)$ for $y\in \O$.  In particular, for $p=2$, the functions $k_{x_1}, \ldots, k_{x_n}$ are the Riesz-representers of $\f_1, \ldots, \f_n$, respectively.  This problem has been addressed in detail in the papers in the papers  \cite{MTN-SP, Nair-ACM,  Nair-discrete} for the case $p=2$.
\eeg 

The following example is an analogue of the above in a more general setting.

\beg 
Let $\X$ and $\Y$ be normed linear spaces, and let $K: \X\to \Y$ be a continuous linear operator.  Given continuous linear functionals $\psi_1, \ldots, \psi_n$ on $\Y$, one may want to determine $f\in \X$ such that 
$$\psi_i(Kf)  =  \xi_i\q\h{for}\q  i=1, \ldots, n.$$
In this case, for each $i\in \{1,\ldots, n\}$, we may  define  $\f_i: \X\to \R$ by 
$$\f_i(f) =  \psi_i(Kf),\q f\in \X.$$
Thus, in this  case, the  training set  is $\{(\f_i, \xi_i): i=1, \ldots, n\}$. 
\eeg

\beg 
Let $\O$  be a  compact subset of $\R^d$ for some $d\in \N$ and for  $1\leq p\leq \infty$, and let  $\X$ be a subspace of $L^p(\O)$.  Suppose we are interested in identifying an $f\in \X$ such that 
 $$\int_{\O_i} f(x)dx =  \xi_i$$
 for some subsets $\O_1, \ldots, \O_n$ of $\O$ and for some real numbers $\xi_1, \ldots, \xi_n$.  In this case, defining 
 $$\f_i(f) = \int_{\O_i} f(x)dx\q\h{for}\q i=1, \ldots, n,$$
 we see, by H\"{o}lder's inequality, that each $\f_i$ is a continuous linear functional on $\X$.  Thus, in the case, the  training set  is $\{(\f_i, \xi_i): i=1, \ldots, n\}$. 
\eeg 
 
\section{The Problem  under  Least-square setting} \label{sec-2} 
 Corresponding to  $\f_1, \ldots, \f_n \in \X^*$, let  $\Phi: \X\to \R^n$ be defined by 
\beq\label{op} 
\Phi(f) = [\f_1(f),  \cdots, \f_n(f)]^T,\q f\in \X.\eeq 
Let  $\xi= (\xi_1, \ldots, \xi_n) \in \R^n$. We consider the operator equation 
\beq\label{op-1} 
\Phi(f) = \xi.
\eeq
Clearly, the above equation has a solution for every choice of $\xi\in \R^n$ if and only if $\Phi$ is onto.   In this context, we observe the following result, which follows using  standard arguments in linear algebra on a  particular type of linear transformation.  Since this result is not available explicitly  in the standard text books on linear algebra, we provide its proof as well.
	
\bt\label{Th-onto-0}
Let $V$ be a vector space over the field $\K$, which is either $\K$ or $\C$,  and let $\psi_1, \ldots, \psi_m$ be linear functionals on $V$. Then the map $F: V\to \K^m$ defined by 
$$F(v) = (\psi_1(v), \ldots, \psi_m(v)),\q v\in V,$$
is onto if and only if $\psi_1, \ldots, \psi_m$ are linearly independent, and in that case,   there exist $v_1, \ldots, v_m$ in $X$ such that 
$$\psi_i(v_j) = \d_{ij},\q i,j\in \{1, \ldots, m\}.$$ 
\et 
	
\bpf 
Taking an orthonormal basis of the range of $\Phi$, \w.r.t. the standard inner product on $\K^m$, and extending it to an orthonormal bases of $\K^n$ (cf. \cite{NS}, Theorem 4.36),  we obtain that  
$F$ is not onto if and only if  there exists a nonzero $n$-tuple  $(\a_1, \ldots, \a_m)  \in \K^m$ such that  $\sum_{i=1}^m \a_i\psi_i(u) = 0$ for all $v\in V$. But, for a  nonzero $(\a_1, \ldots, \a_m)  \in \K^m$,   $\sum_{i=1}^m \a_i\psi_i(v) = 0$ for all $v\in V$  if and only if $\sum_{i=1}^m \a_i\psi_i =0$  if and only if  $\psi_1, \ldots, \psi_m$ are linearly dependent. 
Further,   if $F$ is onto, then  for each $e_j = (\d_{1j}, \ldots, \d_{mj} )\in \K^m$ for $j\in \{1, \ldots, m\}$, there exists $v_j\in V$ such that $F(v_j)=e_j$, equivalently,  $\psi_i(v_j) = \d_{ij}$ for $ i,j\in \{1, \ldots, m\}$. 
\epf

\brem
It is to be remarked  that  the second part of  Theorem \ref{Th-onto-0}   is a  known result (see, e.g., \,\cite{BVL-FA}, Lemma 16.2, Page 283) proved in the context of a \nls $V$. However,  we derived it as a particular case of a general result using simpler arguments in the setting of a vector  space. 
\erem

In view of Theorem \ref{Th-onto-0}, the operator equation (\ref{op-1}) has  a solution for every $\xi\in \R^n$ if and only if $\f_1, \ldots, \f_n$ are linearly independent.  More precisely, we have the following reformulation of  Theorem \ref{Th-onto-0}.

\bt\label{Th-onto}
Given  $\f_1, \ldots, \f_n$ in $\X^*$,  the map $\Phi$  defined in (\ref{op})   is onto if and only if   $\f_1, \ldots, \f_n$ are linearly independent  in $\X^*$, and in that case,  there exists $f_1, \ldots, f_n$ in $\X$ such that 
$$\f_i(f_j) = \d_{ij},\q i,j\in \{1, \ldots, n\}.$$ 
\et

Next, suppose  that $\xi$ is not in the  range  of $\Phi$, so that  equation (\ref{op-1}) has no solution. In this context,  we may look for an $f\in \X$ such that $\Phi(f) \simeq \xi$, that is,  $\Phi(f)$ is an approximation of $\xi$ in some sense.   In this case, we look for a  {\it least-square solution}  for the operator equation (\ref{op-1}),  that is, 
$$\|\Phi(f) - \xi\|_{\R^n} = \inf_{h\in \X}\|\Phi(h)-\xi\|_{\R^n}, $$
where $\|\cdot\|_{\R^n}$ is any norm on $\R^n$ induced by an inner product $\<\cdot, \cdot\>_{\R^n}$. 
For example,  $\<\cdot, \cdot\>_{\R^n}$ can be  the standard inner product on $\R^n$ or a   weighted inner product induced by an $n$-tuple  $(w_1, \ldots, w_n)$ of positive reals, namely, 
\beq\label{ip-n} \<{\bf u, v}\>_{\R^n} = \sum_{i=1}^n w_i{\bf u}_i{\bf v}_i,\q {\bf u, v}\in \R^n.\eeq 
In this case the   corresponding norm $\|\cdot\|_w$ is given by   
$$\|{\bf u}\|_{\R^n} = \Big(\sum_{i=1}^n w_i |{\bf u}_i|^2\Big)^{\frac{1}{2}},\q {\bf u}\in \R^n.$$
In the above, we denoted the $i$-th coordinate of ${\bf u}\in \R^n$ by ${\bf u}_i$. 

We shall show that a least-square solution for (\ref{op-1}) does exist and it is also a solution of an operator equation  of the from 
\beq\label{eq-normal} \Phi^*\Phi(f) = \Phi^*(\xi),\q f\in \X,\, \xi\in \R^n,\eeq 
where  $\Phi^*: \R^n\to \X^*$ is defined by exploiting the inner product structure on $\R^n$, which we shall define in the following theorem (Theorem \ref{Th-adj}).

\brem 
It is to be remarked that, though the notion of  least-square solutions can be defined on  linear spaces with co-domain a Hilbert space, their equivalence to an equation of the form (\ref{eq-normal}) is known only when the the domain space is a Hilbert space (cf. \cite{Nair-linop}, Section 4.2).  Since we are assuming that $\X$ is only a normed linear space, some of the results available for the setting of a Hilbert  space  are to be proved in the setting of a normed linear space. 
\erem 

\brem 
The consideration of a weighted inner product on $\R^n$ found to be useful in the context of certain discrete regularization methods for integral equations of the first kind, where the weights $w_1, \ldots, w_n$ are nothing but the weights used for for a convergent quadrature rule (cf. \cite{Nair-ACM, MTN-SP}).  In the present work, we use such an inner product to get a matrix equation associated with (\ref{eq-normal}) and also to devise a regularization method.  
\erem 

\noi
{\bf{Notation:}} Throughout the paper,  if $X$ and $Y$ are linear spaces and  $T: X\to Y$, is a linear  operator, then  its range  will be denoted by $R(T)$ and its kernel or null-space will be denoted by $N(T)$. Thus, 
$$
R(T)  =  \{Tu: u\in X\},\q N(T) = \{u\in X: Tu=0\}.$$
If $X$ is a  \nlsp, 
then the norm on $X$ will be 
denotd by $\|\cdot\|_X$. If $X$ is an inner product space, then the inner product on $X$  will be denote by  $\<\cdot, \cdot\>_X$ and the induced the norm will be denoted by $\|\cdot\|_X$, and for  $S\subseteq X$,  the orthogonal complement of $S$ will be denoted by $S^\perp$.  Thus, 
$$S^\perp = \{u\in X: \<u, v\>_X=0\h{ for every } v\in S\}.$$
If $X$ is a normed linear space, then the set of all continuous linear functionals on $X$ is denoted by $X^*$, and  for $\f\in X^*$ and $u\in X$, the scalar $\f(u)$ will also be denoted by the {\it duality pair} $\<u, \f\>$.  
Thus,      
$$\<f, \f\> :=  \f(f)\q\forall\, \f\in X^*,\, f\in X.$$
This notation is found to be useful while deducing results associated with the adjoint of an operator from a normed linear space to $\R^n$.
Recall also that $\Phi: \X\to \R^n$ is a continuous linear operator defined by (\ref{op}) \w.r.t. any  inner product $\<\cdot, \cdot\>_{\R^n}$ on $\R^n$, and \w.r.t. this inner product,   we have 
$$\R^n = R(\Phi) \oplus R(\Phi)^\perp.$$ 
The conclusions in the following theorem can be verified easily. 

\bt\label{Th-adj} 
For each ${\bf u}\in \R^n$, the function $f\mapsto   \<\Phi(f), {\bf u}\>_{\R^n}$ belongs to $\X^*$ and 
  $\Phi^*:  \R^n\to \X^*$  defined given by 
$$\<f, \Phi^*({\bf u})\> = \<\Phi(f), {\bf u}\>_{\R^n},\q f\in \X,\, {\bf u}\in \R^n,$$
is a continuous linear operator. 
\et 

We call the operator $\Phi^*,$ obtained  as in the above theorem,  as the {\bf adjoint} of $\Phi$.

\bt\label{adj}
Suppose the  inner product  $\<\cdot, \cdot\>_{\R^n}$ is  the weighted  inner product induced by $(w_1, \ldots, w_n)$. Then $\Phi^*: \R^n\to \X^*$ and $\Phi^*\Phi: \X\to \X^*$ are given by 
$$\Phi^*({\bf u}) = \sum_{j=1}^n w_j {\bf u}_j \f_j,\q {\bf u}\in \R^n,$$ 
$$\Phi^*\Phi(f) =  \sum_{j=1}^n w_j  \<f, \f_j\>  \f_j,$$
respectively.
\et

\bpf 
For $f\in \X$ and ${\bf u}\in \R_w^n$,  we have 
$$
\<f, \Phi^*({\bf u})\> =  \<\Phi(f), {\bf u}\>_{\R^n} 
= \sum_{j=1}^n w_j \< f,   \f_j \> {\bf u}_j = \Big\< f,  \sum_{j=1}^n w_j {\bf u}_j \f_j \Big\>.$$
Thus,  
$$\Phi^*({\bf u}) = \sum_{j=1}^n w_j {\bf u}_j \f_j,\q {\bf u}\in \R^n.$$ 
Therefore, for every  $f\in \X$, we have 
$$\Phi^*\Phi(f) =  \sum_{j=1}^n w_j  \<f, \f_j\>  \f_j.$$
This completes the proof.\epf 

Next we show that (\ref{op-1}) has a least-square solution for every $\xi\in \R^n$. For this, we  shall make use of the following lemma. 

\bl\label{lem-adj}   Let $\Phi^*$ be  the adjoint of $\Phi$, defined as in Theorem \ref{Th-adj}.   Then 
$$ N(\Phi^*) = R(\Phi)^\perp, \q N(\Phi^*\Phi) = N(\Phi).$$
\el 

\bpf
We observe that, for ${\bf u}\in \R_w^n$, 
\beqarray 
\Phi^*({\bf u}) = 0 & \iff &    \<f, \Phi^*({\bf u}) \> = 0\q\forall\, f\in \X,\\
& \iff &    \<\Phi(f),  {\bf u} \>_{\R^n} = 0\q\forall\, f\in \X,\\
& \iff &     {\bf u}\in R(\Phi)^\perp.
\eeqarray 
Hence,  $N(\Phi^*) = R(\Phi)^\perp.$   
Now, to show that $N(\Phi^*\Phi) = N(\Phi)$, first we  observe that 
$N(\Phi) \subseteq N(\Phi^*\Phi)$. To show the other way inclusion, let $f\in N(\Phi^*\Phi)$. Then, we have $\Phi(f) \in N(\Phi^*)=R(\Phi)^\perp$, so that we obtain $\Phi(f) \in  R(\Phi) \cap R(\Phi)^\perp$, resulting in $\Phi(f)=0$, that is, $f\in N(\Phi)$.  Thus, the proof is completed. 
\epf 

Now the theorem on the existence of least-square solution for (\ref{op-1}). 
 
\bt\label{Th-LS}
 There exists $f\in \X$ such that 
\beq\label{eq-LS}\|\Phi(f) - \xi\|_{\R^n} = \inf_{h\in \X}\|\Phi(h)-\xi\|_{\R^n}.\eeq 
Further,  for $f\in \X$, the above equality  holds iff 
$$\Phi^*\Phi(f) = \Phi^*(\xi),$$
where $\Phi^*$ is the adjoint of $\Phi$, defined as in Theorem \ref{Th-adj}.  
\et 

\bpf 
Let $P: \R_w^n\to \R^n$ be the  orthogonal projection onto $R(\Phi)$.  Let  $f\in \X$ be such that  $\Phi(f) = P\xi$. Then, we obtain 
$$\|\xi  - \Phi(f)\|_{\R^n} = \|\xi  - P\xi\|_{\R^n} = 
 \inf_{{\bf u}\in R(\Phi)}\|\xi- {\bf u}\|_{\R^n} =  \inf_{h\in \X}\|\xi-\Phi(h)\|_{\R^n}.$$
Thus,  this $f$ satisfies (\ref{eq-LS}).   
Also, for any $f\in \X$,  using Lemma \ref{lem-adj}, we obtain 
\beqarray \Phi(f) = P\xi  &\iff&  P(\Phi(f)-\xi) = 0 \\
& \iff&  \Phi(f)-\xi\in R(\Phi)^\perp = N(\Phi^*) \\
&\iff & \Phi^*\Phi(f) = \Phi^*(\xi).
\eeqarray 
This completes the proof. 
\epf

\brem
It is to be observed that Theorem \ref{Th-LS}, specially the second part,  is analogous to Theorem 4.5 in \cite{Nair-linop} and Theorem 14.3.2 in \cite{Nair-FA}, which are proved in the case when $\X$ is a Hilbert space. However,  using the inner product structure on $\R^n$, we could  prove  Theorem \ref{Th-LS}  using only the normed linear structure on  $\X$.   
\erem

\section{Identification of a unique least-square solution}

By   Theorem \ref{Th-onto},  if $\f_1, \ldots, \f_n$ are linearly independent, then   the equation   
(\ref{op-1}) 
has a solution in $\X$, but not necessarily unique.    However,  in this case, we can identify a unique solution in a  certain finite dimensional  subspace of $\X$, as shown in the following theorem. 

\bt\label{LS-1} 
Let  $\f_1, \ldots, \f_n$ be  linearly independent  in $\X^*$ and  $f_1, \ldots, f_n$ in $\X$ be   as in Theorem  \ref{Th-onto}. Then  the equation  (\ref{op-1}) has  a unique  solution  $f_0$ belonging to $\h{span}\{f_1, \ldots, f_n\}$, and  it is given by 
$$f_0= \sum_{i=1}^n \xi_if_i.$$
\et 

\bpf 
Let  $f_0 = \sum_{i=1}^n \xi_if_i$. Since $\<f_i, \f_j\> = \d_{ij}$,  we have  
$$\<f_0, \f_j\> =  \sum_{i=1}^n \xi_i \<f_i, \f_j\> = \xi_j,\q j\in \{1, \ldots, n\}. $$
Hence, $\Phi(f_0) = \xi$. To see the uniqueness,  suppose  $f$ is  any arbitrary element in $\h{span}\{f_1, \ldots, f_n\}$ which  satisfies $\Phi(f) = \xi$.  Then $f$ has the representation $f=\sum_{i=1}^n \a_i f_i$ for some $(\a_1, \ldots, \a_n)\in \R^n$,  so that  
$\a_j = \<f, \f_j\>  = \xi_j$ for every $ j\in \{1, \ldots, n\}$. Hence,  we can conclude that $f = f_0$. 
\epf 

Now, suppose  $\Phi$ is not necessarily onto. Let $\psi_1, \ldots, \psi_k$ be linearly independent continuous  linear functionals on $\X$ such that 
\beq\label{span-op} \span\{\psi_1, \ldots, \psi_k\} = \span\{\f_1, \ldots, \f_n\}.\eeq 
By Theorem \ref{Th-onto-0},  there exists $h_1, \ldots, h_k$ in $\X$ such that 
$$\psi_i(h_j) = \d_{ij}\q\h{for}\q i, j\in \{1, \ldots, k\}.$$
Clearly, $h_1, \ldots, h_k$ are linearly independent. Define $Q: \X\to \X$ by 
\beq\label{gen-proj} Q(f) = \sum_{i=1}^k \psi_i(f)h_i,\q f\in \X.\eeq 

\bt\label{projection} 
Let $Q: \X\to \X$  be defined as in (\ref{gen-proj}) and let $\X_0:= \span\{h_1, \ldots, h_k\}$. Then $Q$ is a projection operator with 
$$R(Q)=\X_0,\q N(Q)= N(\Phi),\q\h{and}\q R(\Phi) = \{\Phi(f): f\in \X_0\}.$$
\et 

\bpf Since   $Q(h_j)=h_j$ for every $j\in \{1, \ldots, k\}$, we obtain $Q(f) = f$ for every $f\in \X_0$.  Also,  $f\in \X$,  we have 
\beqarray 
f\in N(Q)  &\iff&  \psi_i(f)=0\q\forall\, i\in \{1, \ldots, k\} \\
&\iff & \f_i(f)=0\q\forall\, i\in \{1, \ldots, n\}\\
&\iff& f\in N(\Phi).\eeqarray
Thus, $Q$ is a projection operator with 
$R(Q)=\X_0$ and $N(Q)= N(\Phi).$

Next, for every $f\in \X$, if $g\in N(\Phi)$ and $h\in \X_0$ are such that $f=g+h$, then we have $\Phi(f) = \Phi(h)$ so that, it follows that $R(\Phi) = \{\Phi(f): f\in \X_0\}$.
\epf 
%\newpage 
Using the above setting, we have the following theorem.

\bt\label{LS-2}
Let $\xi\in \R^n$.  Then,  (\ref{op-1}) has a unique least-square solution in $\X_0$, that is, there  exists a unique    $ f^\dagger \in \X_0$  such that   
$$\Phi^*\Phi(f^\dagger) = \Phi^*(\xi).$$
Further, \w.r.t. $f^\dagger$ is the unique least-square solution  of (\ref{op-1})  such that 
\beq\label{mini-norm} \|f^\dagger \|_0 = \inf\{\|f\|_0: \Phi^*\Phi(f) = \Phi^*(\xi)\},\eeq 
where $\|\cdot\|_0$ is the norm on $\X$, defined by 
$$\|f\|_0 = \|g\|_\X + \|h\|_\X, f\in \X,$$ 
with $(g, h)$  as the unique pair in $N(\Phi)\times \X_0$ such that $f=g+h$.  
\et

\bpf 
Let   $f\in \X$ be any given least-square solution of  the equation  (\ref{op-1}), and let $f^\dagger = Q(f)$ be as defined in (\ref{gen-proj}), where $Q: \X\to \X$ is the projection operator defined as in (\ref{gen-proj}), that is,   
$$f^\dagger = \sum_{i=1}^k \psi_i(f)h_i, $$
where $\psi_1, \ldots, \psi_k$ are as in (\ref{span-op}). 
Clearly, $f^\dagger \in \X_0$. Since $\psi_i(f) = \psi_i(f^\dagger)$ for every $i\in \{1, \ldots,  k\}$, we obtain  $\Phi(f) = \Phi(f^\dagger)$, and hence 
$$\Phi^*\Phi(f^\dagger) = \Phi^*\Phi(f) = \Phi^*(\xi).$$
Now, suppose $\tilde f$ is any least square solution in $\X_0$. Then, we have   
$$\Phi^*\Phi(\tilde f) = \Phi^*(\xi) = \Phi^*\Phi(f^\dagger)$$ 
so that 
$\Phi(\tilde f-f^\dagger) \in N(\Phi^*).$
By Lemma \ref{lem-adj},  $N(\Phi^*) = R(\Phi)^\perp$. Hence, $\Phi(\tilde f-f^\dagger)=0$, that is,   $\tilde f-f^\dagger \in N(\Phi))$.  Since $\tilde f -f^\dagger\in \X_0$, the above equality implies that $\tilde f = f^\dagger$.  This finishes the first part of the theorem.  

For the second part, first  we note that  
$$\|f\|_0 = \|Q(f)\|_\X+\|f-Q(f\|_\X,\q f\in \X.$$
Then, for any least-square solution $f$ of (\ref{op-1}), we have $Q(f) = f^\dagger$ so that 
$$\|f^\dagger\|_0 = \|f^\dagger\|_\X \leq \|f\|_0.$$
Thus, $f^\dagger$ satisfies (\ref{mini-norm}). 
To see the uniqueness of $f^\dagger$ satisfying (\ref{mini-norm}), let $\tilde f$ be also a  least-square solution of (\ref{op-1})  satisfying (\ref{mini-norm}). Then we have 
$$\|f^\dagger\|_\X = \|f^\dagger\|_0 = \|\tilde f\|_0 = \|f^\dagger \|_\X+\|\tilde f-f^\dagger\|_\X.$$
This shows that $\tilde f = f^\dagger$, and the proof is complete. 
\epf

We call the  unique  least-square solution  $f^\dagger$, obtained in the above theorem   as the {\bf generalized solution}  of the equation (\ref{op-1}).

\brem 
It is to be observed that if $\X$ is a Hilbert space, then $f^\dagger$ obtained in Theorem \ref{LS-2} is nothing but $\Phi^\dagger(\xi)$, where $\Phi^\dagger: \R^n\to \X$ is the generalized (Moore-Penrose) inverse of $\Phi$ (cf. \cite{EHN, Nair-linop}).
\erem 

\section{Computation of generalized solution}

In order to consider the problem of computation of $f^\dagger $, which is   obtained as in  Theorem  \ref{LS-2}, we  consider the weighted inner product on $\R^n$ induced by an $n$-tuple $(w_1, \ldots, w_n)$ of positive reals, defined as in (\ref{ip-n}). 
Now, let $h_1, \ldots, h_k$ in $\X$ be as in Theorem  \ref{projection}, and let $f^\dagger\in \X_0:= \h{span}\{h_1, \ldots, h_k\}$   be the unique element in $\X_0$, obtained as in  Theorem  \ref{LS-2},  satisfying 
$$\Phi^*\Phi(f^\dagger ) = \Phi^*(\xi).$$
Let $(\a_1, \ldots, \a_k)\in \R^k$ be the unique element in $\R^k$ such that 
$f^\dagger  = \sum_{j=1}^k \a_jh_j.$
Then using the representations of $\Phi^*$ and $\Phi^*\Phi$, obtained as in Theorem  \ref{adj}, we obtain 
$$\sum_{i=1}^n w_i\<f^\dagger , \f_i\> \f_i =  \sum_{i=1}^n w_i\xi_i  \f_i.  $$
Evaluating the  functions on both sides of the above equation at $h_\ell$ for $\ell\in \{1, \ldots, k\}$, we obtain 
$$\sum_{i=1}^n w_i\<f^\dagger , \f_i\> \<h_\ell, \f_i\> =  \sum_{i=1}^n w_i\xi_i  \<h_\ell, \f_i\>.  $$
Since  $\<f^\dagger , \f_i\>  = \sum_{j=1}^k \a_j\<h_j, \f_i\>$, the above equation takes the form 
$$ \sum_{j=1}^k \Big(\sum_{i=1}^n w_i  \<h_j, \f_i\> \<h_\ell, \f_i\>  \Big)\a_j=  \sum_{i=1}^n w_i\xi_i  \<h_\ell, \f_i\>.$$
Thus, $(\a_1, \ldots, \a_k) \in \R^n$ is obtained by solving the system of equations 
$$\sum_{j=1}^k a_{\ell\,j} \a_j = \b_\ell,\q \ell=1, \ldots, k,$$
where 
\beq\label{matrix-0} a_{\ell\,j} =  \sum_{i=1}^n w_i  \<h_j, \f_i\>  \<h_\ell, \f_i\>,\q  \b_\ell =  \sum_{i=1}^n w_i\xi_i  \<h_\ell, \f_i\>.\eeq 
Thus, the problem of computation of $f^\dagger$ is reduced to the problem of solving the matrix equation 
\beq\label{matrix-1} {\bf A}{\bf u} = {\bf v},\eeq 
for ${\bf u} = (\a_1, \ldots, \a_k) \in \R^n$, where ${\bf A} = (a_{\ell\,j})$ is a $k\times k$ matrix and ${\bf v} = (\b_1, \ldots, \b_k)\in \R^n$.  
\vsqq 
\noi{\bf Algorithm for computation of $f^\dagger $:}
\ben 
\i Solve (\ref{matrix-1}) for ${\bf u} = (\a_1, \ldots, \a_k)$.
\i Write $f^\dagger  = \sum_{j=1}^k \a_j h_j$.
\een 

Before closing this section, let us  observe the following result. 

\bt\label{Th-matrix-pos-self} 
The  $k\times k$ matrix ${\bf A} =(a_{ij})$, where $a_{ij}$ is an in  (\ref{matrix-0}), is symmetric and positive definite \w.r.t. the standard inner product on $\R^k$. Further, $N({\bf A})=0.$
\et

\bpf  Fro the expression for $a_{\ell j}$ as given in  (\ref{matrix-0}), we see that $a_{\ell\,j} = a_{j\ell}$. Thus,  ${\bf   A}$ is symmetric. Now, to see that it is positive definite, let    ${\bf u}\in \R^k$. Then we have 
\beqarray 
\<{\bf A}{\bf u}, {\bf u}\>_{\R^k} &=&     \sum_{\ell=1}^k {\bf  u}_\ell\sum_{j=1}^k   \Big(\sum_{i=1}^n w_i  \<h_j, \f_i\>  \<h_\ell, \f_i\>\Big)      {\bf   u}_j \\
& = & \sum_{i=1}^n  w_i    \Big(\sum_{\ell=1}^k  
\<h_\ell, \f_i\> {\bf  u}_\ell\Big)  \Big( \sum_{j=1}^k      \<h_j, \f_i\>     {\bf   u}_j \Big).
\eeqarray 
Thus, 
\beq\label{pos-def} \<{\bf A}{\bf u}, {\bf u}\>_{\R^k}  = \sum_{i=1}^n  w_i   \Big( \sum_{j=1}^k      \<h_j, \f_i\>     {\bf   u}_j \Big)^2 \geq 0\q\forall\, {\bf u}\in \R^k,
\eeq 
so that   ${\bf A}$ is positive definite.  Also,  in view of  (\ref{pos-def}) and (\ref{span-op}), we have 
\beqarray 
\<{\bf A}{\bf u}, {\bf u}\>_{\R^k}  = 0  
&\iff&  \sum_{j=1}^k      \<h_j, \f_i\>     {\bf   u}_j  = 0 \q\forall\, i=1, \ldots, n,\\
&\iff&  \sum_{j=1}^k      \<h_j, \psi_i\>     {\bf   u}_j  = 0 \q\forall\, i=1, \ldots, k\\
&\iff& {\bf u}_i=0\q\forall\, i=1, \ldots, k \\
&\iff& {\bf u}=0.
\eeqarray
In particular, $N({\bf A})=0$.
\epf

\bcor\label{cor-matrix-pos-self} The equation (\ref{matrix-1})  has a  unique solution ${\bf u}^\dagger$, and there exists a unique ${\bf w}\in \R^k$ such that ${\bf u} ^\dagger = {\bf A}{\bf w}$.
\ecor 

\bpf 
By Theorem \ref{Th-matrix-pos-self}, ${\bf A}: \R^k\to \R^k$  is bijective, and hence the result. 
\epf 

\section{Regularized solution} 

It can happen that the matrix ${\bf A}$ in the equation (\ref{matrix-1}) is ill-conditioned, which means, since ${\bf A}$ is symmetric and positive definite (cf. Theorem \ref{Th-matrix-pos-self}), that the smallest positive eigenvalue of ${\bf A}$ can be too close to $0$. Therefore, it is advisable to consider a regularized approximation of  the solution  ${\bf u}$ of (\ref{matrix-1}), and thereby obtain an approximation for  the generalized solution  $f^\dagger$ of (\ref{op-1}) (cf.  Theorem  \ref{LS-2}). 

In view of Theorem \ref{Th-matrix-pos-self},     the matrix ${\bf A} +\l {\bf I}$ is invertible for every $\l>0$, where ${\bf I}$ is the identity $k\times k$ matrix.  Now,  for $\l>0$,  let ${\bf u}_\l :=(\g_1, \ldots, \g_k)$ be the unique vector in $\R^k$ such that 
\beq\label{matrix-2}({\bf A}+\l {\bf I}) {\bf u}_\l = {\bf v},\eeq 
where ${\bf v}\in \R^k$ is as in  (\ref{matrix-1}).
As a regularized approximation of $f^\dagger$,  we define  
\beq\label{reg-sol} f_\l = \sum_{j=1}^k \g_j h_j\q\h{for}\q \l >0.\eeq

\subsection{Error estimates under noise-free data}

The following theorem not only shows that ${\bf u}_\l$ and $f_\l$ are approximations of ${\bf u}^\dagger $ and  $f^\dagger$, respectively,  for small enough $\l$, but also gives estimates for the errors $\|{\bf u}^\dagger -{\bf u}_\l\|_{\R^k}$ and $\|f^\dagger -f_\l\|_\X$. 

\bt\label{Th-error-noise-free}
Let ${\bf u} ^\dagger = (\a_1, \ldots, \a_k)$ and ${\bf u}_\l = (\g_1, \ldots, \g_k)$ be the solutions of  (\ref{matrix-1}) and (\ref{matrix-2}), respectively,   and  let  $f^\dagger  = \sum_{j=1}^k \a_j h_j$ and $f_\l = \sum_{j=1}^k \g_j h_j$ be 
as   obtained  in Theorem  \ref{LS-2} and in (\ref{reg-sol}), respectively. 
Then    the inequalities 
\beq\label{error-1} \|{\bf u} ^\dagger -{\bf u}_\l\|_{\R^k} \leq \l \|{\bf w}\|_{\R^k},\q \|f^\dagger -f_\l\|_\X  \leq  \Big( \sum_{j=1}^k\|h_j\|_\X^2\Big)^{\frac{1}{2}}\|{\bf u}^\dagger  -{\bf u}_\l\|_{\R^k}\eeq 
hold for every $\l>0$, where ${\bf w}$ is as in Corollary \ref{cor-matrix-pos-self}.
In particular,
\beq\label{error-2}\|{\bf u}^\dagger -{\bf u}_\l\|_{\R^k} =O(\l),   \q   
\|f^\dagger -f_\l\|_\X  = O(\l)  \q\h{as}\q \l \to 0.\eeq 
\et 

\bpf 
Since  ${\bf u}^\dagger $ is the solution of  (\ref{matrix-1}), we have 
$({\bf A} +\l {\bf I}) {\bf u} ^\dagger = \l{\bf u} ^\dagger + {\bf v}.$
This combined with (\ref{matrix-2}) imply that 
$$({\bf A} +\l {\bf I}) ({\bf u}^\dagger  -{\bf u}_\l) = \l{\bf u}^\dagger $$
so that, by the invertibility of ${\bf A} +\l {\bf I}$, we obtain 
$${\bf u} ^\dagger -{\bf u}_\l = \l ({\bf A} +\l {\bf I})^{-1} {\bf u}^\dagger .$$
{By Corollary \ref{cor-matrix-pos-self},  there exists a unique ${\bf w}\in \R^k$ such that   ${\bf u}^\dagger  = {\bf A}{\bf w}$},  so that we obtain 
$$\|{\bf u}^\dagger  -{\bf u}_\l\|_{\R^k} = \l \|({\bf A} +\l {\bf I})^{-1} {\bf A}{\bf w}\|_{\R^k} \leq \l \|{\bf w}\|_{\R^k}.$$
Thus, first inequality in (\ref{error-1}) is proved. For getting the last inequality in the above,  we used the inequality $\|({\bf A} +\l {\bf I})^{-1} {\bf A}\|\leq 1$ (cf. \cite{Nair-linop}). Since 
$$ f^\dagger -f_\l  =    \sum_{j=1}^k (\a_j - \g_j)  h_j,\q  \|{\bf u}^\dagger  - {\bf u}_\l\|_{\R^k}^2 = \sum_{j=1}^k |\a_j - \g_j|^2,$$ 
by Cauchy-Schwarz inequality, we obtain 
$$
\|f^\dagger -f_\l\|_\X  \leq   \sum_{j=1}^k |\a_j - \g_j|\,  \| h_j\|_\X  \leq  \Big( \sum_{j=1}^k  \| h_j\|_\X^2\Big)^{\frac{1}{2}} \|{\bf u} ^\dagger - {\bf u}_\l\|_{\R^k} .$$
Thus, the second inequality in (\ref{error-1}) holds.  Since $\|{\bf u} ^\dagger -{\bf u}_\l\|_{\R^k} \leq \l \|{\bf w}\|_{\R^k}$, we have 
$$
\|f^\dagger -f_\l\|_\X  \leq  \l \|{\bf w}\|_{\R^k} \Big( \sum_{j=1}^k  \| h_j\|_\X^2\Big)^{\frac{1}{2}} .$$
Thus, (\ref{error-2}) is also proved. \epf

\subsection{Error estimates under noisy data}

Next, suppose that the available data  is noisy, that is, we have a noisy data $\tilde \xi$ in place of the actual data $\xi$. Let ${\bf \tilde u}_\l$ be the solution of the equation (\ref{matrix-2}) with ${\bf \tilde v}$ in place of ${\bf v}$, that is, 
\beq\label{matrix-2-noise} ({\bf A}+\l {\bf I}) {\bf \tilde u}_\l = {\bf \tilde v},\eeq 
for   $\l>0$, where 
$${\bf \tilde v} = (\tilde \b_1, \ldots, \tilde \b_k)\in \R^n\q\h{with}\q  \tilde \b_\ell =  \sum_{i=1}^n w_i\tilde \xi_i  \<h_\ell, \f_i\>.$$
Writing 
${\bf \tilde u}_\l :=(\tilde \g_1, \ldots, \tilde \g_k)$ let  
\beq\label{reg-sol-noise} \tilde f_\l = \sum_{j=1}^k \tilde \g_j h_j.\eeq 
Then, we have 
\beq\label{error-noise-1}  \|f_\l-\tilde f_\l\|_\X 
\leq \sum_{j=1}^k|\g_j-  \tilde \g_j|\| h_j\|_\X  \leq  \Big(\sum_{j=1}^k \| h_j\|_\X^2\Big)^{\frac{1}{2}}\|{\bf u}_\l-{\bf \tilde u}_\l\|_{\R^k} .\eeq 

In order to estimate  $\|{\bf u}_\l-{\bf \tilde u}_\l\|_{\R^k}$, first we shall represent the  matrix ${\bf A}$ and  the vector ${\bf v}$ in  an alternative forms.  For this, 
let ${\bf B}  :=(b_{ij})$ be the $n\times k$ matrix with 
$$b_{ij} = \<h_j, \f_i\> \q\h{for}\q i=1, \ldots, n; \, j=1, \ldots, k.$$
Considering ${\bf B}$ as a linear transformation from $\R^k$ to $\R^n$, we have 
$$({\bf  B}{\bf  x})(i) = \sum_{j=1}^k b_{ij}{\bf  x}(j) \q\h{for}\q   {\bf x}\in \R^k;\,  i=1, \ldots, n$$
Here, we denoted the $i$-th coordinate of a vector ${\bf  a}$ by ${\bf  a}(i)$.  Now, let $\R^k$ be endowed with the standard inner product $\<\cdot, \cdot\>_{\R^k}$  and let $\R^n$ be with the weighted inner product $\<\cdot, \cdot\>_w$ induced by $(w_1, \ldots, w_n)$ as defined in (\ref{ip-n}). Then, for ${\bf x}\in \R^k$ and ${\bf y}\in \R^n$,  we have 
\beqarray 
\<{\bf  B}{\bf  x}, {\bf  y}\>_w  &=&  \sum_{i=1}^n  w_i \Big( \sum_{j=1}^k b_{ij}{\bf  x}(j)\Big) {\bf y}(i) =  \sum_{j=1}^k  {\bf  x}(j)  \Big( \sum_{i=1}^n w_i  b_{ij}{\bf y}(i) \Big)  \\
&=& \<{\bf x}, {\bf B}^*{\bf y}\>_{\R^k} ,
\eeqarray 
where  the adjoint operator ${\bf B}^*: \R^n_w\to \R^k$ is given by 
\beq\label{adj-2} ({\bf B}^*{\bf y}) (\ell)=   \sum_{i=1}^n w_i  b_{i\ell}{\bf y}(i) =   \sum_{i=1}^n w_i  \<h_\ell, \f_i\> {\bf y}(i),\eeq 
for $\ell=1, \ldots, k$. Therefore,  ${\bf B}^*{\bf B}: \R^k\to \R^k$ is given by	\beq\label{adj-3}
({\bf B}^*{\bf B} {\bf x})(\ell)  =   \sum_{i=1}^n w_i  b_{i\ell}\Big(  \sum_{j=1}^k b_{ij}{\bf  x}(j) \Big) =\sum_{j=1}^k \Big(   \sum_{i=1}^n w_i  \<h_\ell, \f_i\> \<h_j, \f_i\> \Big){\bf  x}(j)\eeq  
for $\ell=1, \ldots, k$. The representations of  ${\bf B}^*$ and  ${\bf B}^*{\bf B} $ in (\ref{adj-2}) and (\ref{adj-3}), respectively, show that  ${\bf B}^*{\bf B} $ is nothing but the matrix ${\bf A}$ in (\ref{matrix-1}), and the equation (\ref{matrix-1}) gives 
\beq\label{normal-eq-1} 
{\bf B}^*{\bf B} {\bf u}^\dagger  = {\bf B}^* {\bf \xi}.\eeq 

\brem 
We have seen, in Theorem \ref{Th-matrix-pos-self} 
that the  $k\times k$ matrix ${\bf A} =(a_{ij})$  is a symmetric and positive definite matrix. The above discussion  leads to  another proof for this fact.
\erem

From (\ref{matrix-2}), (\ref{matrix-2-noise}) and  (\ref{normal-eq-1}), we have 
\beqarray 
\|{\bf u}_\lambda -{\bf \tilde u}_\l\|_{\R^k}  
&=&  \|({\bf A}+\l {\bf I})^{-1} ({\bf  v} -  {\bf \tilde v})\|_{\R^k} \\
&=&  \|({\bf A}+\l {\bf I})^{-1} {\bf  B}^*(\xi-  \tilde \xi)\|_{\R^k}.
\eeqarray 
But, we know that  (cf. \cite{Nair-linop})
$$\|({\bf A}+\l {\bf I})^{-1} {\bf  B}^* \|= \|({\bf B}^*{\bf B}+\l {\bf I})^{-1} {\bf  B}^* \|\leq \frac{1}{2\sqrt\l}.$$ 
Hence, we obtain 
\beq\label{error-noise-2} \|{\bf u}_\lambda -{\bf \tilde u}_\l\|_{\R^k}   \leq  \frac{\|\xi  -  \tilde \xi\|_w}{2\sqrt\l}\eeq

%%The above discussion together with Theorem \ref{Th-error-noise-free} lead to the  following theorem.

In view of the estimates  in (\ref{error-noise-1}) and   (\ref{error-noise-2}),  and the estimates in Theorem \ref{Th-error-noise-free},  we have the following theorem. 

\bt\label{error-noise}
Let $\tilde \xi$ be a noisy data in place of the actual data $\xi$ and for  $\l>0$, let  ${\bf \tilde u}_\l$ and $\tilde f_\l$ be as in  (\ref{matrix-2-noise}) and (\ref{reg-sol-noise}), respectively.  Then, 
\beq\label{error-3} \|{\bf u}^\dagger -{\bf \tilde u}_\l\|_{\R^k}  \leq  
\l \|{\bf w}\|_{\R^k}   +   \frac{\|\xi  -  \tilde \xi\|_w}{2\sqrt\l} ,\eeq 
\beq\label{error-4} \|f^\dagger -\tilde f_\l\|_\X  \leq \Big(\sum_{j=1}^k \| h_j\|_\X^2\Big)^{\frac{1}{2}}    \|{\bf u}^\dagger -{\bf \tilde u}_\l\|_{\R^k}\eeq 
for all $\l>0$, where ${\bf w}$ is as in Corollary \ref{cor-matrix-pos-self}.  In particular, if  $\|\xi-\tilde \xi\|_w
\leq \d$ for some $\d>0$ and if $\l_{\d}= \d^{2/3}$, then 
\beq\label{error-5} \|{\bf u}^\dagger -{\bf \tilde u}_{\l_\d}\|_{\R^k}  \leq  C_1\d^{2/3},\q  \|f^\dagger -\tilde f_{\l_\d}\|_\X  \leq C_2\d^{2/3},\eeq 
where 
$ C_1 :=   \|{\bf w}\|_{\R^k}   +  1/2$ and  $C_2 :=  C_1   \Big(\sum_{j=1}^k \| h_j\|_\X^2\Big)^{\frac{1}{2}}.$
\et 

\bpf
Since 
$$\|{\bf u} ^\dagger - {\bf \tilde u}_\l\|_{\R^k} \leq  
\|{\bf u} ^\dagger - {\bf u}_\l\|_{\R^k} + \|{\bf u}_\l - {\bf \tilde u}_\l\|_{\R^k}$$
and 
$$\|f^\dagger - \tilde f_\l\|_\X \leq \|f^\dagger - f_\l\|_\X + \|f_\l - \tilde f_\l\|_\X, $$
the inequalities  in (\ref{error-3}) and (\ref{error-4}) follow from the inequalities in  (\ref{error-1}) and (\ref{error-noise-1}). 

{Now, let $\|\xi-\tilde \xi\|_w
\leq \d$ for some $\d>0$, then the  inequalities in (\ref{error-5}) are obtained from (\ref{error-3}) and (\ref{error-4}) by taking $\l={\l_\d}:=\d^{2/3}$.
}
\epf

\brem 
It can be seen easily that $f_\l$ and $\tilde f_\l$ for $\l>0$ defined in (\ref{reg-sol}) and (\ref{reg-sol-noise}), respectively, satisfy the following equations:
$$\Phi^*\Phi(f_\l) + \l  \sum_{j=1}^k \<f_\l, \psi_j\>\psi_j = \Phi^*(\xi),$$
$$\Phi^*\Phi(\tilde f_\l) + \l \sum_{j=1}^k \<\tilde f_\l, \psi_j\>\psi_j = \Phi^*(\tilde \xi).$$
In case the {\nls}  $\X$ is a Hilbert space, then we may  identify $\{\psi_1, \ldots, \psi_k\}$ with an orthonormal basis of $\X_0$ so that the above two equations {\cb are} reduced to  
$$\Phi^*\Phi(f_\l) +  \l f_\l = \Phi^*(\xi)\q\h{and}\q  \Phi^*\Phi(\tilde f_\l) + \l \tilde f_\l = \Phi^*(\tilde \xi),$$
respectively, which are the Tikhonov regularized equations associated with the equation (\ref{op-1}), without noise and with noise in the data, respectively. 
\erem 

\brem 
We observe that the error estimates derived in Theorem \ref{error-noise} is order optimal for  Tikhonov regularization in the setting of Hilbert spaces (cf. \cite{EHN, Nair-linop}). 
\erem 

% {\cred It is to be observed that the problem addressed here  is for a  fixed data space, namely, the space $\R^n$ for a fixed $n\in \N$. Thus, the analysis carried out in this paper does not allow, $n$ to grow indefinitely,  as the the coefficients $C_1$ and $C_2$ in Theorem  \ref{error-noise} depend on $n$. In case the theory has to be applied for  larger data set, one may have to consider a bigger data space $\R^N$ with $N>> n$, and then treat data $(\xi_1, \ldots, \xi_n)\in \R^n$ as $(\xi_1, \ldots, \xi_n, 0, \ldots, 0)\in \R^N$. In any case, for a reasonably small error estimate, $N$ has to be fixed, and cannot be  too large. The investigation of the  effect of $n$ on the constants $C_1$ and $C_2$ in Theorem  \ref{error-noise} is not in the purview of the present work. Such investigation may be possible if the functionals $\f_1, \ldots, \f_n$ arise out of certain specific operators as discussed in the Section 2.} 

\brem  The convergence analysis carried out in Section 6  is  valid when $\l\to 0$ for the noise-free data and when $\d\to 0$ for the noisy data, but not allowing $n$ to be arbitrarily large,   as the  coefficients involved in the estimates can be  very  large.  
The investigation on the  effect of $n$ on the constants   is not in the purview of the present work. Such investigation may be possible under particular  settings in which the functionals $\f_1, \ldots, \f_n$ arise out of certain specific operators as discussed in the Section 2, which can be a topic for a future project.
\erem 

\section{Concluding remarks}

Given a \nls $\X$ and a  training set $\{(\f_i, \xi_i)\in \X^* \times \R: i=1, \ldots, n\}$, we addressed the problem of identifying an  $f\in \X$  such that $\f_i(f)=  \xi_i$  for $i=1, \ldots, n$. This problem is formulated as the problem of solving an operator equation $\Phi(f) = \xi$, where $\Phi: \X\to \R^n$ is defined by $\Phi(f) = (\f_1(f), \ldots, \f_n(f))$  and $\xi = (\xi_1, \ldots, \xi_n)\in \R^n$. In this paper,  this problem is addressed when $\xi \in R(\Phi)$, and also when $\xi\not\in R(\Phi)$ by considering a  least-square solution, and in both the cases, a unique generalized solution $f^\dagger$ is identified in a certain finite dimensional subspace of $\X$. A computational scheme is suggested by formulating the problem in the setting of a matrix equation, which also enabled us  to consider a regularization method to take care of the case when the matrix under consideration is ill-conditioned. Error estimates are derived for the exact data $\xi$ and also when it is noisy, by choosing the regularization parameter appropriately. It is shown that the approach adopted is a generalization to Tikhonov regularization with the domain space a \nlsp. 

%\vsh 
%\noindent{\bf Acknowledgement:} {
%I gratfully acknowledge the anonymous referees for their critical comments and very useful suggestions earlier drafts of the paper, which greatly helped to revise it  thoroughly incorporating many changes and also bringing it into much better form. I also acknowledge the support received from BITS Pilani, K K Birla Goa Campus, where I am a visiting professor after my superannuation from IIT Madras. }

\vsh 
\noindent{\bf Acknowledgement:} {
	I gratfully acknowledge the  support received from BITS Pilani, K K Birla Goa Campus, where I am a visiting professor after my superannuation from IIT Madras. }

\end{document}